\documentclass{article}

\usepackage{amssymb}
\usepackage{amsmath}

\usepackage[dvips]{graphicx}

\begin{document}

 
\begin{center}
\textbf{ Travaux de Husain et al. sur la continuit\'{e} automatique des caract\`{e}res}                                                                                                                                                                                   
\end{center}

\noindent \textbf{}

\begin{center}
M. El Azhari
\end{center}

\noindent \textbf{ } 

\noindent \textbf{Abstract.} We give a survey of Husain, Ng and Liang's results concerning the automatic continuity of algebra homomorphisms. We also give the improvements obtained by Joseph, Akkar, Oudadess, Zelazko and ourselves.

\noindent \textbf{}

\noindent \textit{Mathematics Subject Classification 2010:} 46H40.

\noindent \textit{Key words:} automatic continuity of algebra homomorphisms.
 
\noindent \textbf{} 
 
\noindent \textbf{} 
 
\noindent\textbf{I. Introduction}

\noindent \textbf{} 

\noindent \textbf{} Nous donnons un aper\c{c}u des r\'{e}sultats de T. Husain, S. B. Ng et J. Liang concernant la continuit\'{e} ou la bornitude des homomorphismes d'alg\`{e}bres. Nous donnons aussi les am\'{e}liorations apport\'{e}es par G. A. Joseph, M. Akkar, M. Oudadess, W. Zelazko et nous m\^{e}mes (cf. [2], [3], [15], [18] et [20]).

\noindent \textbf{} Nous obtenons d'abord un th\'{e}or\`{e}me (th\'{e}or\`{e}me III.1.6) de continuit\'{e} s\'{e}quentielle des homomorphismes sur les alg\`{e}bres s\'{e}quentielles. Ce th\'{e}or\`{e}me recouvre et g\'{e}n\'{e}ralise les r\'{e}sultats de T. Husain et S. B. Ng ([12, th\'{e}or\`{e}me 1]), G. A. Joseph    ([15, th\'{e}or\`{e}me 2.1]), M. Oudadess ([18, th\'{e}or\`{e}me 5.2]), et T. Husain ([6, th\'{e}or\`{e}me 2]).   

\noindent \textbf{} Ensuite, Nous obtenons un th\'{e}or\`{e}me (th\'{e}or\`{e}me III.3.3) de continuit\'{e} automatique sur les alg\`{e}bres r\'{e}elles. Comme cons\'{e}quence, on a deux r\'{e}sultats de T. Husain et S. B. Ng ([11, th\'{e}or\`{e}me 1] et [10, th\'{e}or\`{e}me 1]).

\noindent \textbf{} Dans la derni\`{e}re partie consacr\'{e}e aux alg\`{e}bres \`{a} bases, nous donnons deux nouvelles preuves du r\'{e}sultat de T. Husain et J. liang affirmant que si $A$ est une alg\`{e}bre de Fr\'{e}chet \`{a} base orthogonale et inconditionnelle, alors tout caract\`{e}re de $A$ est continu.

\noindent \textbf{}

\noindent \textbf{II. Pr\'{e}liminaires}

\noindent \textbf{}

\noindent \textbf{II.1. Alg\`{e}bres topologiques} 

\noindent \textbf{}

\noindent \textbf{} Une alg\`{e}bre topologique est une alg\`{e}bre (sur $K=R$ ou $C$) munie d'une topologie s\'{e}par\'{e}e, compatible avec sa structure d'espace vectoriel et pour laquelle la multiplication est s\'{e}par\'{e}ment continue. Une alg\`{e}bre topologique est dite localement convexe (en abr\'{e}g\'{e} a.l.c.) si elle munie d'une topologie d'espace localement convexe (en abr\'{e}g\'{e} e.l.c.). Si elle m\'{e}trisable et compl\`{e}te, on dit que c'est une $B_{0}$-alg\`{e}bre. Une alg\`{e}bre topologique $A$ est dite localement multiplicativement convexe (en abr\'{e}g\'{e} a.l.m.c.) si $A$ est munie d'une topologie d\'{e}finie par une famille $(p_{\lambda})_{\lambda\in\Lambda}$ de semi-normes d'espace vectoriel v\'{e}rifiant en outre 
$p_{\lambda}(xy)\leq p_{\lambda}(x)p_{\lambda}(y)$ pour tout $\lambda\in\Lambda$ et tous $x, y$ in $A.$ Une a.l.m.c compl\`{e}te m\'{e}trisable est dite une alg\`{e}bre de Fr\'{e}chet. Soit $A$ une alg\`{e}bre topologique, on note par $M(A)$ l'ensemble des caract\`{e}res continus (non nuls) de $A.$

\noindent \textbf{} Soient $A$ une a.l.c. et $\beta_{A}$ l'ensemble des born\'{e}s, absolument convexes, ferm\'{e}s, idempotents de $A.$ Pour chaque $B\in\beta_{A},$ on note par $A_{B}$ le sous espace vectoriel engendr\'{e} par $B,\,A_{B}$ munie de la jauge $\Vert .\Vert_{B}$ est une alg\`{e}bre norm\'{e}e. Si   $(A_{B},\Vert .\Vert_{B})$ est compl\`{e}te pour tout $B\in\beta_{A},$ on dit que $A$ est pseudocomplete.

\noindent \textbf{}

\noindent \textbf{II.2. Alg\`{e}bres infras\'{e}quentielles et rayon de r\'{e}gularit\'{e}}

\noindent \textbf{}

\noindent \textbf{} Soit $A$ une alg\`{e}bre topologique.

\begin{enumerate}

\item[(i)] $A$ est fortement s\'{e}quentielle s'il existe un voisinage $U$ de $0$ tel que, pour tout $x\in U,\, x^{k}\rightarrow_{k} 0.$

\item[(ii)] $A$ est  s\'{e}quentielle si pour toute suite $(x_{n})_{n}, x_{n}\rightarrow_{n} 0,$ il existe un $x_{m}$ tel que $ x_{m}^{k}\rightarrow_{k} 0. $

\item[(iii)] $A$ est infras\'{e}quentielle si, pour tout born\'{e} $B$ de $A,$ il existe $ \lambda > 0$ tel que, pour tout $x$ de $B$, $(\lambda x)^{k}\rightarrow_{k} 0$    

\end{enumerate}

\noindent \textbf{} On a $(i)\Rightarrow (ii)\Rightarrow (iii).$

\noindent \textbf{•}

\noindent \textbf{•} Soit A une a.l.c. Pour $x\in A,$ on pose:

\noindent \textbf{•} 

\noindent \textbf{•} $\beta(x)=\inf\lbrace r > 0:((r^{-1}x)^{n})_{n}$ est born\'{e}e $\rbrace,$ o\`{u} $\inf\emptyset = \infty$.

\noindent \textbf{•}     
 
\noindent \textbf{•} $\beta(x)$ est la r\'{e}gularit\'{e} de $x$ et $\beta$ le rayon de r\'{e}gularit\'{e}. Un \'{e}l\'{e}ment $x$ est r\'{e}gulier si, et seulement si, $\beta(x) < \infty.$ D'apr\`{e}s [4], on a: 

\begin{enumerate}

\item[(1)] $\beta(\lambda x)=\vert\lambda\vert\beta(x)$ pour tous $x\in A$ and $\lambda\in K$;

\item[(2)] $\infty > \vert\lambda\vert >\beta(x)$ entraine $(\lambda^{-1}x)^{n}\rightarrow_{n} 0.$

\end{enumerate}

\noindent \textbf{•} Soit $A$ une a.l.c. D'apr\`{e}s [18], on a:

\begin{enumerate}

\item[(i)] $A$ est infras\'{e}quentielle si, et seulement si, $\beta$ est born\'{e};

\item[(ii)] $A$ est s\'{e}quentielle si, et seulement si, pour toute suite $(x_{n})_{n}$ convergeant vers $0,$ il existe un certain $m$ tel que $\beta(x_{m}) < 1.$

\item[(iii)] $A$ est fortement s\'{e}quentielle si, et seulement si, $\beta$ est continue.

\end{enumerate}

\noindent \textbf{}

\noindent \textbf{II.3. Alg\`{e}bres topologiques \`{a} bases }

\noindent \textbf{}

\noindent \textbf{} Soit $E$ un espace vectoriel topologique (en abr\'{e}g\'{e} e.v.t.). $E$ est dit \`{a} base s'il existe une suite $(x_{i})_{i\geq 1}$ d'\'{e}l\'{e}ments de $E$ telle que pour tout $x$ de $E$ il existe une suite complexe unique $(\alpha_{i})_{i\geq 1}$ tel que $ x=\lim_{n\rightarrow\infty}\Sigma_{i=1}^{n}\alpha_{i}x_{i}=\Sigma_{i=1}^{\infty}\alpha_{i}x_{i}.$ Pour chaque $i\geq 1$, l'application lin\'{e}aire $ c_{i}:E\rightarrow C,\, c_{i}(x)=\alpha_{i},$ est appel\'{e}e fonction coordonn\'{e}e. $ (x_{i})_{i\geq 1} $ est dite une base inconditionnelle de $E$ si $ \Sigma_{i=1}^{\infty}a_{i}\alpha_{i}x_{i}\in E $ d\`{e}s que $\Sigma_{i=1}^{\infty}\alpha_{i}x_{i}\in E $ et $\vert a_{i}\vert\leq 1. $

\noindent \textbf{} Soit $A$ une alg\`{e}bre topologique \`{a} base $(x_{i})_{i\geq 1},$ on dit que la base  $(x_{i})_{i\geq 1}$ est orthogonale si $x_{i}x_{j}= \delta_{ij}x_{i}$ pour tous $i\geq 1$ et $j\geq 1$, o\`{u} $\delta_{ij}$ est le symbole de Kronecker.

\noindent \textbf{} Soit $E$ un espace de Fr\'{e}chet dont la topologie est d\'{e}finie par la famille\\
 $\lbrace\Vert .\Vert_{r}: r\in N\rbrace$ de seminormes et soit $(x_{n})_{n\geq 1}$ une base inconditionnelle de $E.$ Soient $f\in E^{'}$ (le dual topologique de $E$) et $x\in E,\\
\sup_{f\in \Lambda_{r}}\Sigma_{n=1}^{\infty}\vert c_{n}(x)\vert\vert f(x_{n})\vert <\infty,$ o\`{u} $\Lambda_{r}=\lbrace f\in E^{'}:\sup _{\Vert x\Vert_{r}\leq 1}\vert f(x)\vert\leq 1\rbrace.$\\
Posons $\Vert x\Vert_{r}^{'}=\sup_{f\in \Lambda_{r}}\Sigma_{n=1}^{\infty}\vert c_{n}(x)\vert\vert  
f(x_{n})\vert,$ alors $\lbrace\Vert .\Vert_{r}^{'}: r\in N\rbrace$ est une famille de seminormes d\'{e}finissant la topologie originale de $E.$
 
\noindent \textbf{}

\noindent \textbf{III. Sur les travaux de Husain, Ng et Liang}

\noindent \textbf{}

\noindent \textbf{III.1. Alg\`{e}bres s\'{e}quentielles}

\noindent \textbf{}

\noindent \textbf{} Dans [12], T. Husain et S. B. Ng consid\'{e}rent les alg\`{e}bres s\'{e}quentielles et donnent les r\'{e}sultats suivants:

\noindent \textbf{}

\noindent \textbf{Th\'{e}or\`{e}me III.1.1([12]).} Soit $A$ une a.l.c. s\'{e}quentiellement compl\`{e}te. Si $A$ est s\'{e}quentielle, alors tout caract\`{e}re de $A$ est born\'{e}.
 
\noindent \textbf{}

\noindent \textbf{Corollaire III.1.2([12]).} Soit $A$ une $B_{0}$-alg\`{e}bre s\'{e}quentielle, alors tout caract\`{e}re de $A$ est continu.

\noindent \textbf{}

\noindent \textbf{} La d\'{e}monstration du th\'{e}or\`{e}me III.1.1 est assez longue, elle repose sur le fait suivant: \'{e}tant donn\'{e} $z\in A$ tel que $z^{k}\rightarrow_{k}0,$ on construit une alg\`{e}bre de Banach commutative $(B,\Vert .\Vert)$ telle que $B\subset A,\, z\in B$ et $\Vert z\Vert\leq 1.$
 
\noindent \textbf{}Dans [15], G. A. Joseph donne une d\'{e}monstration courte du th\'{e}or\`{e}me III.1.1 en utilisant quelques notions bornologiques dues \`{a} G. R. Allan. Le r\'{e}sultat de G. A. Joseph est le suivant:

\noindent \textbf{}

\noindent \textbf{Th\'{e}or\`{e}me III.1.3([15]).} Soit $A$ une a.l.c. pseudocompl\`{e}te. Si $A$ est s\'{e}quentielle, alors tout caract\`{e}re de $A$ est born\'{e}.

\noindent \textbf{}

\noindent \textbf{} Le th\'{e}or\`{e}me III.1.3 constitue une am\'{e}lioration du th\'{e}or\`{e}me III.1.1.
 
\noindent \textbf{}

\noindent \textbf{Remarque.} Etant donn\'{e} $z\in A$ tel que $z^{k}\rightarrow_{k}0.$ Soit $C_{z}=\lbrace z,z^{2},...,z^{n},...\rbrace$, $C_{z}$ est un born\'{e} idempotent. On consid\`{e}re $F$ la fermeture de l'enveloppe disqu\'{e}e de $C_{z}$, alors $(A_{F},\Vert .\Vert_{F})$ est une alg\`{e}bre de Banach commutative telle que\\
$A_{F}\subset A,\, z\in A_{F}$ et $\Vert z\Vert_{F}\leq 1,$ ainsi une bonne partie de la d\'{e}monstration de Husain et Ng est simplifi\'{e}e.

\noindent \textbf{}

\noindent \textbf{} Dans [18], M. Oudadess caract\'{e}rise les alg\`{e}bres s\'{e}quentielles \`{a} l'aide du rayon de r\'{e}gularit\'{e}. Comme cons\'{e}quence, il obtient le r\'{e}sultat suivant:

\noindent \textbf{}

\noindent \textbf{Th\'{e}or\`{e}me III.1.4([18]).} Soit $A$ une a.l.c. s\'{e}quentielle et pseudocompl\`{e}te. Alors tout caract\`{e}re de $A$ est s\'{e}quentiellement continu.

\noindent \textbf{}

\noindent \textbf{} Le th\'{e}or\`{e}me III.1.4 recouvre et g\'{e}n\'{e}ralise les th\'{e}or\`{e}mes III.1.1 et III.1.3. Dans [6], T. Husain donne le r\'{e}sultat suivant:

\noindent \textbf{}

\noindent \textbf{Th\'{e}or\`{e}me III.1.5([6]).} Soient $A, B$ deux a.l.c. s\'{e}quentiellement compl\`{e}tes, $A$ \'{e}tant s\'{e}quentielle et $B$ satisfait \`{a} la condition suivante:\\
(C): Pour toute suite $(y_{n})_{n}$ dans $B,\,y_{n}\neq 0,\,y_{n}\nrightarrow 0,$ il existe une suite $(f_{m})_{m}$ de caract\`{e}res de $B$ tel que $\inf _{n,m}\vert f_{m}(y_{n})\vert =\epsilon > 0.$ \\
Alors tout homomorphisme $T:A\rightarrow B$ est s\'{e}quentiellement continu.

\noindent \textbf{}

\noindent \textbf{} Le th\'{e}or\`{e}me III.1.4 ne peut pas se d\'{e}duire du th\'{e}or\`{e}me III.1.5 car le corps des complexes $C$ ne satisfait pas \`{a} la condition (C), en effet soit $(y_{n})_{n}$ la suite de $C$ d\'{e}finie par $y_{n}=\frac{1}{n}$ si $n$ est pair et $y_{n}=1$ si $n$ est impair. Il est clair que $y_{n}\neq 0$ et  $y_{n}\nrightarrow 0.$ Pour toute suite de caract\`{e}res de $C,$ on a                  $\inf _{n,m}\vert f_{m}(y_{n})\vert = 0$ car pour tout $m,\,f_{m}$ est ou bien l'application nulle ou bien l'application identique.
 
\noindent \textbf{} En modifiant la condition (C) par la condition (D), on obtient un th\'{e}or\`{e}me qui recouvre et g\'{e}n\'{e}ralise les th\'{e}or\`{e}mes III.1.1, III.1.3, III.1.4 et III.1.5.

\noindent \textbf{}

\noindent \textbf{Th\'{e}or\`{e}me III.1.6.} Soient $A, B$ deux a.l.c., $A$ \'{e}tant s\'{e}quentielle et pseudocomplete et $B$ satisfait \`{a} la condition suivante:\\
(D): Pour toute suite $(y_{n})_{n}$ dans $B,\,y_{n}\neq 0,\,y_{n}\nrightarrow 0,$ il existe un caract\`{e}re $f$ de $B$ tel que $ f(y_{n})\nrightarrow 0. $\\
Alors tout homomorphisme $T:A\rightarrow B$ est s\'{e}quentiellement continu.

\noindent \textbf{}

\noindent \textbf{Preuve.} Supposons que $T$ n'est pas s\'{e}quentiellement continu. Alors il existe une suite $(x_{n})_{n}$ dans $A,\,x_{n}\rightarrow 0$ mais $Tx_{n}\nrightarrow 0.$ Posons $y_{n}=\frac{1}{2}Tx_{n},$ on peut supposer que $y_{n}\neq 0$ pour tout $n.$ D'apr\`{e}s (D), il existe un caract\`{e}re $f$ de $B$ tel que $f(y_{n})\nrightarrow 0.$ On peut trouver une sous suite $(y_{n_{k}})_{k}$ de $(y_{n})_{n}$ telle que $\inf _{k}\vert f(y_{n_{k}})\vert = \epsilon > 0.$ Posons $z_{n_{k}}=\epsilon^{-1}x_{n_{k}},\,z_{n_{k}}\rightarrow_{k} 0.$ Comme $A$ est s\'{e}quentielle, il existe $z\in\lbrace z_{n_{k}}\rbrace_{k},\,z=z_{n_{k_{0}}}$ pour un certain $k_{0},$ tel que $z^{i}\rightarrow_{i} 0.$ Il est clair que $f\circ T$ est un caract\`{e}re de $A.$ Soit $F$ la fermeture de l'enveloppe disqu\'{e}e de $\lbrace z,z^{2},...,z^{n},...\rbrace.\,(A_{F},\Vert .\Vert_{F})$ est une alg\`{e}bre de Banach commutative telle que $A_{F}\subset A,\,z\in A_{F}$ et $\Vert z\Vert_{F}\leq 1.$
On a $\vert f\circ T(x)\vert\leq\Vert x\Vert_{F}$ pour tout $x\in A_{F},$ en particulier                           $\vert f\circ T(z)\vert\leq\Vert z\Vert_{F}\leq 1.$ Par construction,\\
$\vert f\circ T(z)\vert =\vert f(\epsilon^{-1}Tx_{n_{k_{0}}})\vert = 2\epsilon^{-1}\vert f(y_{n_{k_{0}}})\vert\geq 2,$ ce qui est impossible.

\noindent \textbf{}

\noindent \textbf{Remarque.} Toute alg\`{e}bre v\'{e}rifiant la condition (C) satisfait \`{a} la condition (D) et de m\^{e}me le corps des complexes satisfait \`{a} la condition (D).

\noindent \textbf{}

\noindent \textbf{III.2. Alg\`{e}bres infras\'{e}quentielles}

\noindent \textbf{}

\noindent \textbf{} Dans [6], T. Husain consid\`{e}re les alg\`{e}bres infras\'{e}quentielles et donne le r\'{e}sultat suivant:

\noindent \textbf{}

\noindent \textbf{Th\'{e}or\`{e}me III.2.1([6]).} Soit $A$ une a.l.c. s\'{e}quentiellement compl\`{e}te. Si $A$ est infras\'{e}quentielle, alors tout caract\`{e}re de $A$ est born\'{e}.

\noindent \textbf{}

\noindent \textbf{} La d\'{e}monstration du th\'{e}or\`{e}me III.2.1 repose sur le fait d\'{e}j\`{a} cit\'{e}: \'{e}tant donn\'{e} $z\in A$ tel que $z^{k}\rightarrow_{k}0,$ on construit une alg\`{e}bre de Banach commutative $(B,\Vert .\Vert)$ telle que $B\subset A,\, z\in B$ et $\Vert z\Vert\leq 1.$

\noindent \textbf{} Dans [18], M. Oudadess caract\'{e}rise les alg\`{e}bres infras\'{e}quentielles \`{a} l'aide du rayon de r\'{e}gularit\'{e} $\beta.$ Comme cons\'{e}quence, il obtient le r\'{e}sultat suivant:

\noindent \textbf{}

\noindent \textbf{Th\'{e}or\`{e}me III.2.2([18]).} Soit $A$ une a.l.c. $\beta$-r\'{e}guli\`{e}re pseudocomplete. Alors l'ensemble $M^{\ast}$ des caract\`{e}res de $A$ est equiborn\'{e}, i.e., pour tout born\'{e} $B$ de $A,\,\bigcup\lbrace f(B): f\in M^{\ast}\rbrace$ est born\'{e}. En particulier, tout caract\`{e}re de $A$ est born\'{e}.

\noindent \textbf{}

\noindent \textbf{} Le th\'{e}or\`{e}me III.2.2 g\'{e}n\'{e}ralise le th\'{e}or\`{e}me III.2.1.

\noindent \textbf{}

\noindent \textbf{III.3. Alg\`{e}bres r\'{e}elles}

\noindent \textbf{}

\noindent \textbf{} L'outil fondamental dans cette partie est la technique de Do. Sin. Sya ([19]) et sa g\'{e}n\'{e}ralisation. Dans [17], S. B. Ng et S. Warner g\'{e}n\'{e}ralisent la technique de Do. Sin. Sya de la fa\c{c}on suivante:

\noindent \textbf{}

\noindent \textbf{Th\'{e}or\`{e}me III.3.1([17]).} Soit $K$ un corps archim\'{e}dien. Soit $(H,+)$ un groupe topologique ab\'{e}lien, m\'{e}trisable et complet. Soit $s:H\rightarrow H$ une application continue telle que $s(0)=0.$ Si $f$ est un homomorphisme de $(H,+)$ dans $(K,+)$ tel qu'il existe $v\in N,\,f(x)^{2}\leq v.f(s(x))$ pour tout $x$ de $H,$ alors $f$ est continu.

\noindent \textbf{}

\noindent \textbf{} Dans [17], l'\'{e}nonc\'{e} du th\'{e}or\`{e}me III.3.1 est donn\'{e} lorsque $K=R,$ mais le r\'{e}sultat reste vrai lorsque on remplace $R$ par un corps archim\'{e}dien. Comme cons\'{e}quence du th\'{e}or\`{e}me III.3.1, on a:

\noindent \textbf{}

\noindent \textbf{Th\'{e}or\`{e}me III.3.2(Lemme de Do. Sin. Sya).} Soit $A$ une alg\`{e}bre topologique m\'{e}trisable compl\`{e}te. Soit $B$ un $R$-sous espace vectoriel de $A.$ On suppose que $B$ est ferm\'{e} et $\lbrace x^{2}: x\in B\rbrace\subset B.$ Soit $f$ une forme lin\'{e}aire sur $B$ telle que   $f(x)^{2}\leq f(x^{2})$ pour tout $x$ de $B,$ alors $f$ est continue.

\noindent \textbf{}

\noindent \textbf{} Nous donnons un th\'{e}or\`{e}me dont la d\'{e}monstration repose essentiellement sur la g\'{e}n\'{e}ralisation de la technique de Do. Sin. Sya. Comme cons\'{e}quence, nous obtenons deux 
th\'{e}or\`{e}mes de T. Husain et S. B. Ng (th\'{e}or\`{e}me 1 de [11] et th\'{e}or\`{e}me 1 de [10]).

\noindent \textbf{}

\noindent \textbf{Th\'{e}or\`{e}me III.3.3.} $A$ et $B$ \'{e}tant deux $R$-e.v.t. m\'{e}trisables complets. Soient $s:A\rightarrow A,\,h:B\rightarrow B$ deux applications, $s$ \'{e}tant continue et $s(0)= 0.$
Consid\'{e}rons $I_{h}=\lbrace g\in B^{\ast}: g(h(x))=g(x)^{2}$ pour tout $x$ de $B\rbrace$ ($B^{\ast}$ dual alg\'{e}brique de $B$) et supposons que $I_{h}$ est non vide. Si $B$ satisfait \`{a} la condition:\\
(P) pour toute suite $(y_{n})_{n}$ dans $B,\,y_{n}\neq 0,\,y_{n}\nrightarrow 0,$ il existe $f\in I_{h}$ tel que $f(y_{n})\nrightarrow 0;$\\
alors toute application lin\'{e}aire $T$ de $A$ dans $B$ telle que $T(s(x))=h(Tx)$ pour tout $x$ de $A,$ est continue.

\noindent \textbf{}

\noindent \textbf{Preuve.} Supposons que $T$ n'est pas continue. Il existe une suite $(x_{n})_{n}$ dans   $A,\,x_{n}\rightarrow 0$ mais $Tx_{n}\nrightarrow 0.$ On peut supposer que $x_{n}\neq 0$ et $Tx_{n}\neq 0$
pour tout $n\geq 1.$ Par hypoth\`{e}se, il existe $f\in I_{h}$ tel que $f(Tx_{n})\nrightarrow 0.$ On peut construire \`{a} partir de la suite $(x_{n})_{n},$ une suite $(a_{m})_{m}$ telle que $a_{m}\rightarrow 0$   
et $\inf _{m}f(Ta_{m})=\epsilon > 0.$ On a $\epsilon^{-1}a_{m}\rightarrow 0,$ ainsi $s(\epsilon^{-1}a_{m})\rightarrow 0$ car $s$ est continue. Posons $y_{m}=s(\epsilon^{-1}a_{m})$ pour tout $m\geq 1.$ Alors $f(Ty_{m})= f(T(s(\epsilon^{-1}a_{m})))= f(h(T(\epsilon^{-1}a_{m})))= f(T(\epsilon^{-1}a_{m}))^{2}\geq 1$ pour tout $m\geq 1.$ \\
On d\'{e}finit $g_{k}:A^{k+1}\rightarrow A$ pour $k\geq 0:$\\
$g_{0}(y_{1})= y_{1}$\\
$g_{1}(y_{1},y_{2})= y_{1}+s(y_{2})$\\
$\cdots$ \\
$g_{k}(y_{1},...,y_{k+1})= g_{1}(y_{1},g_{k-1}(y_{2},...,y_{k+1})).$\\
En utilisant la m\^{e}me construction faite dans [17], on peut d\'{e}finir une sous suite $(z_{k})_{k\geq 0}$ de $(y_{m})_{m\geq 1}$ telle que pour tout $k\geq 0,\,(g_{p-k}(z_{k},...,z_{p}))_{p\geq k}$ est une suite de Cauchy. Soit $c_{k}= \lim _{p\rightarrow \infty}g_{p-k}(z_{k},...,z_{p}).$ Par d\'{e}finition\\
$g_{p-k}(z_{k},...,z_{p})= z_{k}+s(g_{p-(k+1)}(z_{k+1},...,z_{p})),$ on obtient $c_{k}= z_{k}+s(c_{k+1}),$\\
ainsi $Tc_{k}= Tz_{k}+T(s(c_{k+1})).$ On a $f(Tc_{k})= f(Tz_{k})+f(T(s(c_{k+1})))= f(Tz_{k})+f(h(Tc_{k+1}))= f(Tz_{k})+f(Tc_{k+1})^{2}\geq 1+f(Tc_{k+1})^{2}$ pour tout $k\geq 0.$ D'ou $f(Tc_{0})\geq 1+f(Tc_{1})^{2}\geq 2+f(Tc_{2})^{2}\geq k+f(Tc_{k})^{2},$ i.e. $f(Tc_{0})\geq k$ pour tout $k\geq 0,$ ce qui est absurde.

\noindent \textbf{}

\noindent \textbf{} Comme cons\'{e}quence, on a:

\noindent \textbf{}

\noindent \textbf{Th\'{e}or\`{e}me III.3.4.} Soient $A, B$ deux $R$-alg\`{e}bres topologiques m\'{e}trisables et compl\`{e}tes. On suppose que $B$ satisfait \`{a} la condition:\\
(D): Pour toute suite $(y_{n})_{n}$ dans $B,\,y_{n}\neq 0,\,y_{n}\nrightarrow 0,$ il existe un caract\`{e}re $f$ de $B$ tel que $ f(y_{n})\nrightarrow 0.$\\
Alors toute application lin\'{e}aire $T:A\rightarrow B$ v\'{e}rifiant $T(x^{2})= (Tx)^{2}$ pour tout $x$
de $A,$ est continue.

\noindent \textbf{}

\noindent \textbf{Preuve.} On consid\`{e}re $s:A\rightarrow A,\,s(x)= x^{2}$ et $h:B\rightarrow B,\,h(x)= x^{2}.$ Remarquons que $f$ de la condition (D) est dans $I_{h}.$ On applique alors le th\'{e}or\`{e}me III.3.3.

\noindent \textbf{}

\noindent \textbf{Th\'{e}or\`{e}me III.3.5.} Soit $A$ une $R$-a.l.m.c. s\'{e}quentiellement compl\`{e}te. Soit $B$ une $R$-alg\`{e}bre topologique m\'{e}trisable compl\`{e}te satisfaisant \`{a} la condition (D) du th\'{e}or\`{e}me III.3.4. Alors toute application lin\'{e}aire $T$ de $A$ dans $B$ v\'{e}rifiant   $T(x^{2})= (Tx)^{2},\,x\in A,$ est born\'{e}e.

\noindent \textbf{}

\noindent \textbf{Preuve.} On applique le th\'{e}or\`{e}me III.3.3 et le th\'{e}or\`{e}me de structure de M. Akkar ([1]) affirmant que si $A$ est une a.l.m.c. s\'{e}quentiellement compl\`{e}te, alors $A$ est bornologiquement limite inductive d'alg\`{e}bres de Fr\'{e}chet.

\noindent \textbf{}

\noindent \textbf{Remarques.}\begin{enumerate}
\item Les th\'{e}or\`{e}mes III.3.4 et III.3.5 sont des am\'{e}liorations des th\'{e}or\`{e}mes 1 de [11] et 1 de [10].
\item Dans l'\'{e}nonc\'{e} du th\'{e}or\`{e}me III.3.3, on peut remplacer $R$ par un corps archim\'{e}dien.
\item $C$ est une $R$-alg\`{e}bre de Banach, mais $C$ ne satisfait pas \`{a} la condition (D) car le seul caract\`{e}re r\'{e}el de $C$ est l'application nulle de $C$ dans $R.$
\item On peut remplacer $B$ dans les th\'{e}or\`{e}mes III.3.4 et III.3.5 par $R;$ c'est une $R$-alg\`{e}bre de Banach qui satisfait \`{a} la condition (D).
\end{enumerate}

\noindent \textbf{}

\noindent \textbf{} Comme application du th\'{e}or\`{e}me III.3.3, on a:

\noindent \textbf{}

\noindent \textbf{Corllaire III.3.6.} Soit $A$ un $R$-e.v.t. m\'{e}trisable et complet. Soit $s:A\rightarrow A$ une application continue avec $s(0)= 0$. Consid\'{e}rons $I_{s}=\lbrace g\in A^{\ast}: g(s(x))=g(x)^{2}$ pour tout $x$ de $A\rbrace$ et supposons que $I_{s}$ est non vide. Si $A$ satisfait \`{a} la condition:\\
(P): Pour toute suite $(y_{n})_{n}$ dans $A,\,y_{n}\neq 0,\,y_{n}\nrightarrow 0,$ il existe $f\in I_{s}$ tel que $ f(y_{n})\nrightarrow 0;$\\
alors $A$ poss\'{e}de une unique topologie d'e.v.t. m\'{e}trisable complet.

\noindent \textbf{}

\noindent \textbf{Corollaire III.3.7.} Soit $A$ une $R$-alg\`{e}bre topologique m\'{e}trisable compl\`{e}te satisfaisant \`{a} la condition:\\
(D): Pour toute suite $(y_{n})_{n}$ dans $A,\,y_{n}\neq 0,\,y_{n}\nrightarrow 0,$ il existe un caract\`{e}re $f$ de $A$ tel que $ f(y_{n})\nrightarrow 0.$\\
Alors $A$ poss\'{e}de une unique topologie d'alg\`{e}bre m\'{e}trisable compl\`{e}te.

\noindent \textbf{}

\noindent \textbf{Preuve.} On consid\`{e}re l'application $s:A\rightarrow A,\,s(x)= x^{2}$ et on applique le corollaire III.3.6.

\noindent \textbf{}

\noindent \textbf{Remarque.} Le corollaire III.3.7 est une am\'{e}lioration du corollaire 2 de [11].

\noindent \textbf{}

\noindent \textbf{III.4. Alg\`{e}bres \`{a} bases}

\noindent \textbf{}

\noindent \textbf{} Dans [8] et [9], T. Husain et J. Liang consid\'{e}rent les alg\`{e}bres \`{a} bases
et donnent les r\'{e}sultats suivants:

\noindent \textbf{}

\noindent \textbf{Th\'{e}or\`{e}me III.4.1([9]).} Soit $(A,(p_{i})_{i\geq 1})$ une alg\`{e}bre de Fr\'{e}chet \`{a} base $(x_{i})_{i\geq 1}$ v\'{e}rifiant:\\
(1) $x_{i}x_{j}= x_{j}x_{i}= x_{j},\,i\leq j; $\\
(2) $p_{i}(x_{i})\neq 0,\,p_{i}(x_{i+1})= 0.$\\
Alors tout caract\`{e}re de $A$ est continu.

\noindent \textbf{}

\noindent \textbf{Th\'{e}or\`{e}me III.4.2([8]).} Soit $A$ une alg\`{e}bre de Fr\'{e}chet \`{a} base orthogonale et inconditionnelle. Alors tout caract\`{e}re de $A$ est continu.

\noindent \textbf{}

\noindent \textbf{} Les d\'{e}monstrations de ces deux th\'{e}or\`{e}mes sont longues, calculatoires et n'utilisent aucun r\'{e}sultat de la th\'{e}orie spectrale des a.l.m.c. Dans [3], nous avons donn\'{e} des d\'{e}monstrations courtes des deux th\'{e}or\`{e}mes III.4.1 et III.4.2 en utilisant une caract\'{e}risation du spectre ponctuel ([16, corollaire 5.6]) dans les a.l.m.c. commutatives et compl\`{e}tes.

\noindent \textbf{}

\noindent \textbf{Remarque.} Dans [2], on a montr\'{e} que toute alg\`{e}bre \`{a} base $(x_{i})_{i\geq 1}  $ v\'{e}rifiant les conditions (1) et (2) du th\'{e}or\`{e}me III.4.1 est isomorphe alg\'{e}briquement et topologiquement \`{a} l'alg\`{e}bre  $C^{N}.$

\noindent \textbf{}

\noindent \textbf{} Nous donnons deux nouvelles preuves, \'{e}galement courtes, du   
Th\'{e}or\`{e}me III.4.2 en utilisant le lemme de Do. Sin. Sya (th\'{e}or\`{e}me III.3.2).

\noindent \textbf{}

\noindent \textbf{1\`{e}re preuve du Th\'{e}or\`{e}me III.4.2.} On pose $B=\lbrace x=\Sigma_{k=1}^{\infty}\alpha_{k}x_{k}\in A: \alpha_{k}\in R\,$ pour tout $k\geq 1\rbrace = \bigcap_{n=1}^{\infty}c_{n}^{-1}(R).\,B$ v\'{e}rifie les conditions du th\'{e}or\`{e}me III.3.2. Soit $f$
un caract\`{e}re de $A.$ S'il existe $k\geq 1$ tel que $f(x_{k})\neq 0,$ alors $f$ est continu ([8]). Supposons que $f(x_{k})= 0$ pour tout $k\geq 1.$ Soit $x=\Sigma_{k=1}^{\infty}\alpha_{k}x_{k}\in A.$ On consid\`{e}re les suites complexes $(a_{k})_{k\geq 1}$ et $(b_{k})_{k\geq 1}$ d\'{e}finies par:\\
$a_{k}= 0$ si $\alpha_{k}= 0,\,a_{k}=\alpha_{k}^{-1}Re(\alpha_{k})$ si $\alpha_{k}\neq 0;$\\
$b_{k}= 0$ si $\alpha_{k}= 0,\,b_{k}=\alpha_{k}^{-1}Im(\alpha_{k})$ si $\alpha_{k}\neq 0.$\\
$\vert a_{k}\vert\leq 1$ et $\vert b_{k}\vert\leq 1$ entrainent que $x_{1}=\Sigma_{k=1}^{\infty}a_{k}\alpha_{k}x_{k}\in A$ et $x_{2}=\Sigma_{k=1}^{\infty}b_{k}\alpha_{k}x_{k}\in A.$ On remarque que $x_{1}$ et $x_{2}$ sont des \'{e}l\'{e}ments de $B.$ On a $x= x_{1}+x_{2}$ et $f(x)= f(x_{1})+f(x_{2}).$ Le th\'{e}or\`{e}me III.3.2 entraine que la restriction de $f$ \`{a} $B$ est nulle, d'ou $f$ est nul.

\noindent \textbf{}

\noindent \textbf{} La deuxi\`{e}me preuve du th\'{e}or\`{e}me III.4.2 consiste \`{a} montrer que tout alg\`{e}bre $A,$ de Fr\'{e}chet \`{a} base orthogonale et inconditionnelle, peut \^{e}tre munie d'une involution continue et que tout caract\`{e}re de $A$ est hermitien. Avant de donner la deuxi\`{e}me preuve, nous avons les r\'{e}sultats suivants:

\noindent \textbf{}

\noindent \textbf{Proposition III.4.3.} Soit $A$ une a.l.m.c. compl\`{e}te involutive \`{a} base orthogonale et hermitienne, alors tout caract\`{e}re de $A$ est hermitien.

\noindent \textbf{}

\noindent \textbf{Preuve.} D'apr\`{e}s la proposition 8 de [7], chaque fonction coordonn\'{e}e est hermitienne. Soit $H$ le sous espace vectoriel r\'{e}el de $A$ form\'{e} par les \'{e}l\'{e}ments hermitiens de $A.$ On a $H=\lbrace x\in A: x^{\ast}= x\rbrace =\lbrace x\in A: c_{n}(x^{\ast})= c_{n}(x)$  
pour tout $n\geq 1\rbrace = \lbrace x\in A: \overline{c_{n}(x)}= c_{n}(x)$  
pour tout $n\geq 1\rbrace = \lbrace x\in A: c_{n}(x)\in R\,$   
pour tout $n\geq 1\rbrace.$ Soit $f$ un caract\`{e}re de $A$ et $x\in H.$ Si $f(x)\neq 0,$ il existe $c_{k}\in\lbrace c_{n}: n\geq 1\rbrace$ ([13]) tel que $f(x)= c_{k}(x)\in R.$ Il s'ensuit que la restriction de $f$ \`{a} $H$ est \`{a} valeurs r\'{e}elles. Soit $x\in A,$ en utilisant le fait que $A= H+iH,$ on obtient que $f(x^{\ast})= \overline{f(x)}.$

\noindent \textbf{}

\noindent \textbf{Th\'{e}or\`{e}me III.4.4.} Soit $A$ une alg\`{e}bre de Fr\'{e}chet involutive \`{a} base orthogonale et hermitienne, alors tout caract\`{e}re de $A$ est continu.

\noindent \textbf{}

\noindent \textbf{Preuve.} Soit $f$ un caract\`{e}re de $A.$ D'apr\`{e}s la proposition III.4.3, $f$ est hermitien. $f$ se prolonge \`{a} $A^{+}$ (alg\`{e}bre obtenue par adjonction de l'unit\'{e}) en une forme positive, alors $f$ est continu (th\'{e}or\`{e}me 4 de [17]).

\noindent \textbf{}

\noindent \textbf{2\`{e}me preuve du Th\'{e}or\`{e}me III.4.2.} Si $x= \Sigma_{n=1}^{\infty}\alpha_{n}x_{n}\in A,$ alors $\bar{x}= \Sigma_{n=1}^{\infty}\bar{\alpha}_{n}x_{n}$
converge dans $A,$ ceci provient du fait que la base $(x_{n})_{n\geq 1}$ est inconditionnelle ([14]). On consid\`{e}re l'application $v:A\rightarrow A,\,v(x)= \bar{x},\,v$ est une involution sur $A.$ On remarque que $\Vert x\Vert_{r}^{'}=\Vert\bar{x} \Vert_{r}^{'}$ pour tous $r\in N$ et $x\in A,$ ainsi $v$ est continue. On a $\bar{x}_{n}= x_{n}$ pour tout $n\geq 1,$ d'ou la base $(x_{n})_{n\geq 1}$ est hermitienne. Il s'ensuit d'apr\`{e}s le th\'{e}or\`{e}me III.4.4 que tout caract\`{e}re de $A$ est continu.

\noindent \textbf{}

\noindent \textbf{} Dans [20], en utilisant un r\'{e}sultat de R. Arens ([5, th\'{e}or\`{e}me 6.3]), W. Zelazko obtient le th\'{e}or\`{e}me suivant:

\noindent \textbf{}

\noindent \textbf{Th\'{e}or\`{e}me III.4.5([20]).} Soit $A$ une alg\`{e}bre de Fr\'{e}chet commutative. Si $M(A)$ est auplus d\'{e}nombrable, alors tout caract\`{e}re de $A$ est continu.

\noindent \textbf{}

\noindent \textbf{} Comme cons\'{e}quence, on a:

\noindent \textbf{}

\noindent \textbf{Corollaire III.4.6([20]).} Soit $A$ une alg\`{e}bre de Fr\'{e}chet \`{a} base orthogonale, alors tout caract\`{e}re de $A$ est continu.

\noindent \textbf{}

\noindent \textbf{Remarque.} Le corollaire III.4.6 est une bonne am\'{e}lioration du th\'{e}or\`{e}me III.4.2 car la base de $A$ n'est pas n\'{e}cessairement inconditionnelle.

\noindent \textbf{}

\noindent \textbf{} L'auteur remercie Messieurs les Professeurs M. Akkar et M. Oudadess dont les remarques et suggestions ont contribu\'{e} \`{a} l'am\'{e}lioration du manuscrit.

\noindent \textbf{}   

\noindent \textbf{}  
   
\noindent \textbf{R\'{e}f\'{e}rences}

\noindent \textbf{}

\noindent [1] M. Akkar, Sur la structure des alg\`{e}bres topologiques localement multiplicativement convexes, C. R. Acad. Sc. Paris, 279 (1974), S\'{e}rie A, 941--944.
 
\noindent [2] M. Akkar, M. El Azhari et M. Oudadess, Th\'{e}or\`{e}mes de structure sur certaines alg\`{e}bres de Fr\'{e}chet, Ann. Sc. Math. Quebec, 11 (1987), 245--252.

\noindent [3] M. Akkar, M. El Azhari et M. Oudadess, Continuit\'{e} des caract\`{e}res dans les alg\`{e}bres de Fr\'{e}chet \`{a} bases, Canad. Math. Bull., 31 (1988), 168--174.
  
\noindent [4] G. R. Allan, A spectral theory for locally convex algebras, Proc. London Math. Soc., 15 (1965), 399--421.

\noindent [5] R. Arens, Dense inverse limit rings, Michigan Math. J., 5 (1958), 169--182.

\noindent [6] T. Husain, Infrasequential topological algebras, Canad. Math. Bull., 22 (1979), 413--418.

\noindent [7] T. Husain, Positive functionals on topological algebras with bases, Math. Japonica, 28 (1983), 683--687.

\noindent [8] T. Husain and J. Liang, Multiplicative functionals on Fr\'{e}chet algebras with bases, Canad. J. Math., 29 (1977), 270--276.

\noindent [9] T. Husain and J. Liang, Continuity of multiplicative linear functionals on Fr\'{e}chet algebras with bases, Bull. Soc. Roy. Sc. Li\`{e}ge, 46 (1977), 8--11. 

\noindent [10] T. Husain and S. B. Ng, Boundedness of multiplicative linear functionals, Canad. Math. Bull., 17 (1974), 213--215.

\noindent [11] T. Husain and S. B. Ng, On continuity of algebra homomorphism and uniqueness of metric topology, Math. Zeit., 139 (1974), 1--4.

\noindent [12] T. Husain and S. B. Ng, On the boundedness of multiplicative and positive functionals, J. Austral. Math. Soc.,21 (1976), 498--503.
 
\noindent [13] T. Husain and S. Watson, Topological algebras with orthogonal Schauder bases, Pacific J. Math., 91 (1980), 339--347.

\noindent [14] T. Husain and S. Watson, Algebras with unconditional orthogonal bases, Proc. Amer. Math. Soc., 79 (1980), 539--545.

\noindent [15] G. A. Joseph, Multiplicative functionals and a class of topological algebras, Bull. Austral. Math. Soc., 17 (1977), 391--399.

\noindent [16] E. A. Michael, Locally multiplicatively convex topological algebras, Mem. Amer. Math. Soc., 11 (1952).

\noindent [17] S. B. Ng and S. Warner, Continuity of positive and multiplicative functionals, Duke Math. J., 39 (1972), 281--284.

\noindent [18] M. Oudadess, Rayon de r\'{e}gularit\'{e} dans les alg\`{e}bres infras\'{e}quentielles, Canad. J. Math., 36 (1984), 84--94.

\noindent [19] Do. Sin. Sya, On semi normed rings with an involution, Izv. Akad. Nauk. SSSR, 23 (1959), 509--528.

\noindent [20] W. Zelazko, Functional continuity of commutative m-convex $B_{0}$-algebras with countable maximal ideal spaces, Colloq. Math., 51 (1987), 395--399.

\noindent \textbf{}  

\noindent \textbf{} 

\noindent \textbf{} Ecole Normale Sup\'{e}rieure

\noindent \textbf{} Avenue Oued Akreuch

\noindent \textbf{} Takaddoum, BP 5118, Rabat

\noindent \textbf{} Morocco
 
\noindent \textbf{} 

\noindent \textbf{} E-mail:  mohammed.elazhari@yahoo.fr

\end{document}